\documentclass[amsthm]{elsart}
\usepackage{yjsco}
\usepackage{natbib}
\usepackage{amsmath, amssymb, diagrams}

\theoremstyle{plain}
\newtheorem{lemma}[thm]{Lemma}
\newtheorem{question}{Question}
\newtheorem{remark}[thm]{Remark}

\theoremstyle{definition}
\newtheorem*{example*}{Example}

\newcommand{\pd}[1]{\operatorname{pd}\!\left(#1\right)}
\newcommand{\codim}[1]{\operatorname{ht}\!\left(#1\right)}
\newcommand{\grade}[1]{\operatorname{grade}(#1)}
\newcommand{\ass}[1]{\operatorname{Ass}\!\left(#1\right)}

\newcommand{\type}[2]{\ensuremath{\langle e=#1 \,|\, \lambda=#2\rangle}}

\newcommand{\cm}{Cohen-Macaulay}
\newcommand{\unm}{^\mathrm{unm}}
\newcommand{\m}{\mathfrak{m}}
\newcommand{\p}[1]{\mathbb{P}^{#1}}
\newcommand{\fgh}{\ensuremath{(f,g,h)}}
\newcommand{\pdfgh}{\pd{R/\fgh}}
\def\sfrac#1#2{\hbox{$^{#1} \!/\! _{#2}$}}
\renewcommand{\leq}{\leqslant}
\renewcommand{\geq}{\geqslant}

\date{}

\begin{document}
\begin{frontmatter}

\title{Bound on the projective dimension\\of three cubics}

\author{Bahman Engheta}
\address{Wilshire Fixed Income Analytics, \\
a division of Wilshire Associates Incorporated, \\
Santa Monica, CA 90401, U.S.A.}
\ead{engheta@gmail.com}

\begin{abstract}
We show that given any polynomial ring $R$ over a field and any ideal $J \subset R$ which is generated by three cubic forms, the projective dimension of $R/J$ is at most 36. 
We also settle the question whether ideals generated by three cubic forms can have projective dimension greater than four, by constructing one with projective dimension equal to five. 
\end{abstract}


\end{frontmatter}

\section{Introduction}	\label{section:intro}

Throughout this paper, unless stated otherwise, $R$ denotes \emph{any} polynomial ring over an arbitrary field $k$, say $R=k[X_1,\ldots,X_n]$ where $n$ is not specified, and all ideals are homogeneous. 
Consider the following question posed by Michael E.~Stillman. 

\begin{question}[Stillman {\cite[Problem 3.14]{ps}}]
\label{q}
Is there a bound, independent of $n$, on the projective dimension of ideals in $R=k[X_1,\ldots,X_n]$ which are generated by $N$ homogeneous polynomials of given degrees $d_1,\ldots,d_N$?
\end{question}

Unlike the Hilbert Syzygy Theorem which bounds the projective dimension of an ideal by the dimension $n$ of the underlying ring $R$, this question concerns the existence of a uniform bound on the projective dimension of $R/J$ where neither the ring $R$ nor the ideal $J \subset R$ are fixed, but merely the number of generators of $J$ and the degrees of those generators. 
Equivalently, the above question could be phrased as 
\begin{multline*}
	\sup_n \big\{ \, \pd{R/J} ~~ | ~~ J \subset R=k[X_1,\ldots,X_n] \mbox{ is an ideal} \\
	\mbox{generated by } N \mbox{ forms of degrees } d_1, \ldots, d_N \, \big\} < \infty \;?
\end{multline*}
where $\pd{R/J}$ denotes the projective dimension of $R/J$ over $R$. 

Recall that a rather straightforward construction of Burch \cite{b} shows how already three-generated ideals can have arbitrary large projective dimension. 
Burch's construction, however, comes at the cost of increasing degrees of the generators. 
The assumptions on the number of generators and their degrees are thus easily seen to be necessary. 

Question \ref{q} is further motivated by the notable fact that it is equivalent to the very same question posed about the Castelnuovo-Mumford regularity of ideals in polynomial rings: 
\emph{Is there a bound on the regularity of an ideal solely in terms of the number of its generators and the degrees of those generators?} 
See \cite[Section 1.3]{thesis} for a proof of this equivalence following an argument due to Caviglia. 

In this paper we consider the case $N=3$, $d_1=d_2=d_3=3$, and show that if $J$ is an ideal generated by three cubic forms, then $\pd{R/J} \leq 36$. 
Our goal is to establish the existence of such a bound and not necessarily to obtain the best bound possible; in all likelihood, the bound of 36 is far from being sharp. 
In fact, until recently there were no known examples of three cubics with projective dimension greater than 4. 
In Section~\ref{section:pdfgh=5} we exhibit the only construction known to date which yields three cubics whose projective dimension equals 5. 

The approach presented here is informed by previous work \cite{e} of the author, 
wherein connections to the unmixed part of $I$ and to ideals linked to the unmixed part of $I$ were established -- see Theorems \ref{thm:general}, \ref{thm:unmix-quad}, and \ref{thm:link-quad}.

\subsection{Preliminaries}

{\bf Notation.} 
We denote by $\mathfrak m$ the homogeneous maximal ideal $(X_1,\ldots,X_n)$ of $R$. 
For an ideal $J$, $\codim{J}$ denotes the height of $J$ and $J\unm$ the unmixed part of $J$, that is, the intersection of those primary components $Q$ of $J$ with $\codim{Q} = \codim{J}$. 
By $\lambda(R/J)$ we denote the length of $R/J$ and by $e(R/J)$ its multiplicity at $\mathfrak m$. 
One has $e(R/J) = e(R/J\unm)$ and the associativity formula for multiplicities: 
\begin{equation}		\label{eqn:assoc}
	e(R/J) ~ = \!\!\!\!\!\!
	\sum_{\substack{P \, \in \, \mathrm{Spec}(R) \\ \dim(R/P) \, = \, \dim(R/J) }} 
	\!\!\!\!\!\!
	e(R/P) \; \lambda(R_P/J_P).
\end{equation}
With the associativity formula~\eqref{eqn:assoc} in mind, we adopt the following notation in order to easily refer to an ideal with given multiplicity and number of primary components of minimal height: 
We say that an ideal $J$ is of type 
\[
	\type{a_1,\ldots,a_m}{b_1,\ldots,b_m} 
\]
if $J$ has exactly $m$ associated primes of minimal height with multiplicities $a_1,\ldots,a_m$ and locally at each of those primes $R/J$ has length $b_1,\ldots,b_m$, respectively. 
So $R/J$ has multiplicity $\sum_{i=1}^m a_i\,b_i$ by \eqref{eqn:assoc}. 
Note that an ideal and its unmixed part are of the same type and there are only finitely many possible types for an unmixed ideal of fixed multiplicity. 
For example, prime ideals are of type $\type a1$ and primary ideals are of type $\type ab$. 

\medskip%

The following proposition classifies all height two unmixed ideals of multiplicity two. 
Of interest to us are those of type \type12 which are described in part~(iv). 

\begin{prop}[Engheta {\cite[Proposition~11]{e}}]		\label{prop:class}
Let $R$ be a polynomial ring over a field and let $I \subset R$ be a homogeneous height two unmixed ideal of multiplicity two. 
Then $\pd{R/I} \leq 3$ and $I$ is one of the following ideals. 
\begin{enumerate}
	\item[(i)] A prime ideal generated by a linear form and an irreducible quadric. 
	\item[(ii)] $(x,y) \cap (x,v)=(x,yv)$ with independent linear forms $x,y,v$. 
	\item[(iii)] $(x,y) \cap (u,v)=(xu,xv,yu,yv)$ with independent linear forms $x,y,u,v$. 
	\item[(iv)] The $(x,y)$-primary ideal $(x,y)^2 + (ax+by)$ with independent linear forms $x,y$ and forms $a,b \in \m$ such that $x,y,a,b$ form a regular sequence. 
	\item[(iv$^\circ$)] $(x,y^2)$ with independent linear forms $x,y$. 
\end{enumerate}
\end{prop}

One of the key results in \cite{e} stated that if $J \subset R$ is a three-generated ideal of height two and $I' \subset R$ an ideal linked to the unmixed part of $J$, then $\pd{R/J} \leq \pd{R/I'}+1$. 
We generalize this fact in Theorem~\ref{thm:general} and give a simpler proof. 
To this end, we will need the following elementary lemma. 

\begin{lemma}	\label{lemma:colon-unmix}
If $K$ is an unmixed ideal, then $K:J = K:J\unm$ for any ideal $J$ with $\codim{J} \geq \codim{K}$. 
\end{lemma}
\begin{pf}
As $K:J\unm \subseteq K:J$, it suffices to check the claim locally at every $P \in \ass{R/(K:J\unm)}$. 
As $K$ is unmixed, $\ass{R/(K:J\unm)} \subseteq \ass{R/K}$ and $\codim{P} = \codim{K}$. 
By our assumption, $\codim{P} \leq \codim{J}$ and the claim follows from $J\unm_{_P} = J^{}_{_P}$. 
\end{pf}

\begin{thm}		\label{thm:general}
Let $R$ be a regular local ring and let $J$ be an $N$-generated ideal of $R$ of height $N-1$. 
If $\,\mathbf z = z_{_1}, \ldots, z_{_{N-1}}$ is a maximal regular sequence in $J$, then 
\[
	\pd{R/J} \leq \pd{ \sfrac{R}{(\mathbf z):J} } + 1
\]
and equality holds if and only if $R/J$ is not \cm{}, that is, if and only if $\, \pd{R/J} \geq N$. 
\end{thm}
\begin{pf}
Let $J=(f_1,\ldots,f_N)$ with $\codim{f_1,\ldots,f_{N-1}}=N-1$ and let $\mathbf z$ be a maximal regular sequence in $J$. 
By Lemma~\ref{lemma:colon-unmix}, $(\mathbf z):J = (\mathbf z):J\unm$, that is, $(\mathbf z):J$ is linked to the unmixed part of $J$. 
As any two links of an ideal in a Gorenstein ring have the same (finite or infinite) projective dimension, we have $\pd{ \sfrac{R}{(\mathbf z):J\unm} } = \pd{ \sfrac{R}{(f_1,\ldots,f_{N-1}):J\unm} }$. 
So it suffices to prove the claim for $\mathbf z = f_1,\ldots,f_{N-1}$. 

Notice that $(f_1,\ldots,f_{N-1}):J = (f_1,\ldots,f_{N-1}):f_N$. 
This yields the short exact sequence 
\[
	0 \longrightarrow \frac{R}{(f_1,\ldots,f_{N-1}):J} \xrightarrow{\cdot~f_{_N}} \frac{R}{(f_1,\ldots,f_{N-1})} \longrightarrow \frac{R}{J} \longrightarrow 0, 
\]
of which the middle term $R/(f_1,\ldots,f_{N-1})$ is minimally resolved by the Koszul complex on the elements $f_1,\ldots,f_{N-1}$ and has projective dimension $N-1$. 
Since one has $\pd{ \sfrac{R}{(f_1,\ldots,f_{N-1}):J} } \geq \grade{(f_1,\ldots,f_{N-1}):J} = N-1$, it follows that $\pd{R/J} \leq \pd{ \sfrac{R}{(f_1,\ldots,f_{N-1}):J} } + 1$, as claimed. 

If $R/J$ is not \cm{}, then $\pd{R/J} \geq N$ and we also have the reverse inequality $\pd{ \sfrac{R}{(f_1,\ldots,f_{N-1}):J} } \leq \pd{R/J} - 1$. 
And if $R/J$ is \cm{}, then $J$ is unmixed and $(\mathbf z):J$ is linked to $J$. 
In particular, $\sfrac{R}{(\mathbf z):J}$ is \cm{} as well and $\pd{R/J} = \pd{ \sfrac{R}{(\mathbf z):J} }$. 
\end{pf}

We recall the following theorem which allows us to focus our attention on those ideals whose unmixed part is generated in degree 3 or higher. 

\begin{thm}[Engheta {\cite[Theorem~16]{e}}]		\label{thm:unmix-quad}
Let $R$ be a polynomial ring over a field and let $J \subset R$ be an ideal generated by three cubics. 
If the unmixed part of $J$ contains a quadric form, then $\pd{R/J} \leq 4$. 
\end{thm}

\section{The projective dimension of three cubics}

Let $f,g,h \in R$ be three cubic forms. 
In this section we prove that the projective dimension of $R/\fgh$ is bounded above by 36. 
$I = \fgh\unm$ will denote the unmixed part of the ideal $\fgh$ and $I'$ will be used to denote an ideal which is linked to $I$. 

By \cite[Remark~2]{e} we may assume that $\fgh$ has height two. 
And clearly, we may assume that $f,g,h$ are minimal generators. 
This in turn implies that the multiplicity $e(R/\fgh)$ is at most 8 --- cf.~\cite[Lemma~8]{e}. 
It was shown in \cite{e} that $\pdfgh \leq 3$ if $e(R/\fgh) = 1$, $\pdfgh \leq 4$ if $e(R/\fgh) = 2$, and $\pdfgh \leq 16$ if $e(R/\fgh) = 3$. 

If $e(R/\fgh) = 7$, then we let $p_1,p_2$ be two cubics in $I=\fgh\unm$ which form a regular sequence and we consider the link $I'=(p_1,p_2):I$ which has multiplicity $9-7=2$. 
By Proposition~\ref{prop:class} we have $\pd{R/I'} \leq 3$ and it follows from Theorem~\ref{thm:general} that $\pdfgh \leq 4$. 
Similarly, if $e(R/\fgh) = 8$, then the link $I'$ has multiplicity one and thus $R/I'$ is \cm, that is, $\pd{R/I'}=2$ and $\pdfgh \leq 3$ by Theorem~\ref{thm:general}. 

There remain the cases of multiplicity four, five, and six which will require most of our attention. 
In the following theorem we summarize our results. 

\begin{thm}
If $f,g,h$ are three cubic forms in a polynomial ring $R$ over a field, then $\pdfgh \leq 36$. 
More precisely, with $I=\fgh\unm$, 
\begin{enumerate}
	\item[(a)] If $\,\codim{f,g,h}=3$, or if $\,\codim{f,g,h}=2$ and $I$ contains a linear form, 
	\\%
	then $\pdfgh \leq 3$. 
	\quad%
	(See~\cite[Proposition~6]{e}.) 
	
	\item[(b)] If $\,\codim{f,g,h}=1$ or if $I$ contains a quadric, then $\pdfgh \leq 4$. 
	\\%
	(See~Theorem~\ref{thm:unmix-quad}.) 
	
	\item[(c)] Suppose $e(R/\fgh) \leq 5$ and let $I'$ be an ideal which is linked to $I$ via a complete intersection generated by cubics. 
	If $I'$ contains a quadric, then $\pdfgh \leq 4$. 
	\quad%
	(See~Theorem~\ref{thm:link-quad}.) 
	
	\item[(d)] Below we give a breakdown of the established bounds by multiplicity. 
	\[ \begin{array}{r|r}
		\mathrm{multiplicity~of} & \mathrm{bound~on}
		\\
		$R/\fgh$ & \pdfgh
		\\ \hline
		1, 8 & 3
		\\
		2, 7 & 4
		\\
		3 & 16
		\\
		4 & 36
		\\
		5, 6 & 20
	\end{array} \]
\end{enumerate}
\end{thm}

\subsection{Multiplicity four}
\label{subsection:mult4}

For the case of multiplicity four, we prove Proposition~\ref{prop:bound-36} which supplies a bound of 36 for $\pdfgh$ whenever the ideal $\fgh$ has multiplicity $\geq 2$ along a codimension two linear subspace, that is, whenever $R/I$ has length $\geq 2$ locally at an associated prime of multiplicity one. 
To this end, we will need the following lemma. 

\begin{lemma}	\label{lemma:3-quadrics}
Three quadrics which minimally generate an ideal of height $\leq 2$ can be expressed entirely in terms of 8 linear forms, unless two of the quadrics share a common linear factor. 
\end{lemma}
\begin{pf}
Let $q_1,q_2,q_3$ be three quadrics. 
The statement is evident if $(q_1, q_2, q_3)$ has height one. 
If $\codim{q_1, q_2, q_3} = 2$, then it is easily seen that $(q_1, q_2, q_3)$ has multiplicity $\leq 3$ --- cf.~\cite[Lemma~8]{e}. 
We pass to the unmixed part of $(q_1,q_2,q_3)$ and consider each case separately. 

Let $I$ denote the unmixed part of the ideal $(q_1, q_2, q_3)$ and note that $\codim I = 2$. 
If $e(R/I)=1$, then $I$ is generated by two independent linear forms $x,y$ and $q_i=l_{i1}x+l_{i2}y$ with $i=1,2,3$ and linear forms $l_{i1},l_{i2}$. 
So $q_1, q_2, q_3$ can be expressed in terms of 8 linear forms $l_{i1},l_{i2},x,y$. 

If $e(R/I)=2$, then, by Proposition~\ref{prop:class}, $I$ is one of the following ideals: 

\noindent (i)~ 
$I=(x,q)$ with a linear form $x$ and an irreducible quadric $q$. 
Then $q_i=l_i x+\alpha_i q$ with linear forms $l_i$ and field coefficients $\alpha_i$ for $i=1,2,3$. 
As $\codim{q_1, q_2, q_3} = 2$, the coefficients $\alpha_i$ must not be all zero; 
say $\alpha_3 \neq 0$. 
Replacing $q_1$ by $q_1-\frac{\alpha_1}{\alpha_3} q_3 = (l_1-\frac{\alpha_1}{\alpha_3}l_3)x$ and $q_2$ by $q_2-\frac{\alpha_2}{\alpha_3} q_3 = (l_2-\frac{\alpha_2}{\alpha_3}l_3)x$, they both become divisible by the linear form $x$ and we are done.

\noindent (ii)~ 
$I=(x,yv)$ with independent linear forms $x,y,v$. 
Then $q_i=l_i x+\alpha_i yv$ with linear forms $l_i$ and field coefficients $\alpha_i$ for $i=1,2,3$. 
So $q_1, q_2, q_3 \in k[l_1,l_2,l_3,x,y,v]$. 

\noindent (iii)~ 
$I=(xu,xv,yu,yv)$ with independent linear forms $x,y,u,v$. 
Clearly, we have $q_1, q_2, q_3 \in k[x,y,u,v]$. 

\noindent (iv)~ 
$I=(x,y)^2 + (ax+by)$ with independent linear forms $x,y$ and elements $a,b \in \m$ such that $x,y,a,b$ form a regular sequence. 
As $I$ is the unmixed part of an ideal generated by quadrics, we must have $\deg(ax+by)=2$, for otherwise $I=(x,y)^2$. 
So, $a$ and $b$ are linear and $q_1, q_2, q_3 \in k[a,b,x,y]$. 

\noindent (iv$^\circ$)~ 
$I=(x,y^2)$ with independent linear forms $x,y$. 
In analogy to part~(ii) above, $q_1, q_2, q_3 \in k[l_1,l_2,l_3,x,y]$. 

It remains the case $e(R/I)=3$. 
By the associativity formula~\eqref{eqn:assoc} there are five cases to consider. 
(These cases were discussed in detail in \cite[Section~4]{e}.) 
In three of those cases, $I$ is contained in an ideal generated by two linear forms and, as argued above in the case of multiplicity one, the quadrics $q_1,q_2,q_3$ can be expressed in terms of 8 linear forms. 
We consider the remaining two cases: 

$I$ is a homogeneous prime ideal of minimal multiplicity. As such, $I$ is generated by the $2 \times 2$ minors of a $3 \times 2$ matrix of indeterminates --- cf.~\cite{eh}. 
That is, $I$ is generated by three quadrics in at most six variables, and therefore the same holds for $(q_1,q_2,q_3)$. 

$I$ is primary to $(x,y)$ with independent linear forms $x,y$ and $\lambda\!\left( (R/I)_{(x,y)} \right) = 3$. 
Either $I=(x,y)^2$ or $I$ is generated by $(x,y)^3$ plus additional terms of the form $c_jx+d_jy$ with $(c_j,d_j) \not\subset (x,y)^2$. 
In the former case we are done, as $q_1, q_2, q_3 \in k[x,y]$. 
In the latter case we first rule out the possibility that one of the terms $c_jx+d_jy$ may be linear: if so, then $I=(x,y^3)$ after a linear change of coordinates and thus $(q_1, q_2, q_3) \subset (x)$, a contradiction, since $\codim{q_1, q_2, q_3}=2$. 

So now we have $(q_1, q_2, q_3) \subseteq (c_jx+d_jy)$ with $\deg(c_jx+d_jy) \geq 2$. 
Write $q_i = \sum_j \alpha_{ij} (c_jx+d_jy)$ with field coefficients $\alpha_{ij}$ where $\alpha_{ij}=0$ whenever $\deg(c_jx+d_jy)>2$. 
Then $l_{i1} := \sum_j \alpha_{ij} c_j$ and $l_{i2} := \sum_j \alpha_{ij} d_j$ are linear and $q_1, q_2, q_3$ can be expressed in terms of 8 linear forms $l_{i1},l_{i2},x,y$. 
\end{pf}

\begin{prop}		\label{prop:bound-36}
Let $f,g,h$ be three cubic forms which minimally generate an ideal of height two. 
Suppose that $\fgh$ has a component primary to an ideal $P=(x,y)$ with independent linear forms $x,y$ and $\lambda\!\left( \left( \sfrac{R}{\fgh} \right)_{\!P} \right) \geq 2$. 
Then $\pdfgh \leq 36$. 
\end{prop}
(In our notation, the hypothesis of the proposition simply states that if $\fgh$ is of type \type{a_1,\ldots,a_m}{b_1,\ldots,b_m}, then $a_i=1$ and $b_i \geq 2$ for some $i$.) 
\begin{pf}
Let $Q$ denote the $P$-primary component of $\fgh$, that is, $\fgh \subseteq Q \subsetneq P$ and $\fgh_{_P}=Q_{_P} \subsetneq P_{_P}$. 
We have $e(R/Q) = \lambda(R_{_P}/Q_{_P}) \geq 2$. 
If $Q \subseteq P^2$, then the cubics $f,g,h$ can be expressed in terms of the quadrics $x^2,xy,y^2$ using no more than 9 linear forms $l_i$, in which case $f,g,h \in k[x,y,l_i]$ and $\pdfgh \leq 11$. 
So we may assume that $Q$ contains additional terms of the form $cx+dy$ where $(c,d) \not\subset P$. 
Consequently, the Hilbert function of $(R/Q)_{_P}$ is given by $(\underbrace{1,1,1,\ldots,1}_{e(R/Q)\mathrm{~times}})$. 
(We caution that in addition to $P^{^{e(R/Q)}}$ and the above mentioned terms of the form $cx+dy$, the ideal $Q$ may contain other terms as minimal generators --- cf.~example in \cite[Section~3]{e}.)

Now consider the ideal $I := Q:P^{^{e(R/Q)-2}}$ whose Hilbert function, locally at $P$, is given by $(1,1)$. 
That is, $I$ is a $P$-primary ideal of multiplicity two. 
By parts~(iv) and (iv$^\circ$) of Proposition~\ref{prop:class}, $I = P^2 + (ax+by)$ with elements $a,b$ such that $\codim{x,y,a,b}>3$. 
(The term $ax+by$ need not be the same as the term $cx+dy$ above.) 
In other words, either $x,y,a,b$ form a regular sequence or $(a,b)$ is the unit ideal, in which case we may take $I$ to be $(x,y^2)$. 

Note that $\fgh \subseteq Q \subseteq P^2 + (ax+by)$. 
In what follows, we exploit this inclusion to place $f,g,h$ inside a subring of $R$ generated by a bounded number of linear forms (or by a regular sequence), which will in turn give a bound for $\pdfgh$. 

If $\deg(ax+by)=4$, then $\fgh \subseteq P^2$ and $\pdfgh \leq 11$ as shown above. 
(Strictly speaking, this case is ruled out by our assumption that $Q \neq P^2$.) 

If $\deg(ax+by)=3$, then we may assume without loss of generality that $h=ax+by$ and $f,g \in P^2$.
Indeed, as $\fgh \subseteq P^2 + (ax+by)$, there are nine linear forms $l_{ij}$ and field coefficients $\alpha, \beta, \gamma$ such that 
\[
	\begin{pmatrix}f \\ g \\ h\end{pmatrix}
	=
	\begin{pmatrix}
		l_{11} & l_{12} & l_{13} & \alpha \\
		l_{21} & l_{22} & l_{23} & \beta \\
		l_{31} & l_{32} & l_{33} & \gamma \\
	\end{pmatrix}
	\begin{pmatrix}x^2 \\ xy \\ y^2 \\ ax+by\end{pmatrix}.
\]
If $\alpha=\beta=\gamma=0$, then $\fgh \subseteq P^2$ and we are done; 
so we may assume $\gamma \neq 0$.
Replacing $f$ by $f-\frac{\alpha}{\gamma}h$ and $g$ by $g-\frac{\beta}{\gamma}h$, we have $f,g \in P^2$. 
And relabeling $(l_{31}x+l_{32}y+\gamma a)$ as $a$ and $(l_{33}y+\gamma b)$ as $b$, we can write $h=ax+by$ where $x,y,a,b$ still form a regular sequence. 

Setting $L:=(l_{11},l_{12},l_{13},l_{21},l_{22},l_{23})$, we consider the following two cases: 
If $a$ and $b$ share a common factor modulo $L+P$, then $\pdfgh \leq 27$.
Indeed, if $a \equiv a'c$ and $b \equiv b'c$ modulo $L+P$ with linear forms $a',b',c$, then $a-a'c$ can be written in terms of $x,y,l_{11},\ldots,l_{23}$ using eight linear forms $u_1,\ldots,u_8$ and the same holds for $b-b'c$ with eight linear forms $v_1,\ldots,v_8$. 
Thus, the cubics $f,g,h$ are in the subring $k[x,y,l_{11},\ldots,l_{23},a',b',c,u_1,\ldots,u_8,v_1,\ldots,v_8]$ and $\pdfgh \leq 27$. 
If on the other hand $a$ and $b$ do not have a common factor modulo $L+P$, then they form a regular sequence modulo $L+P$.
That is, the generators of $L+P$ along with $a,b$ form a regular sequence of length at most 10 and $\pdfgh \leq 10$.

If $\deg(ax+by)=2$, then the cubics $f,g,h$ can be expressed in terms of the quadrics $x^2,xy,y^2,ax+by$ using no more than 12 linear forms $l_{ij}$. 
So $f,g,h \in k[x,y,a,b,l_{ij}]$ and $\pdfgh \leq 16$. 

It remains the case where $I=(x,y^2)$. 
Here we have three linear forms $l_1, l_2, l_3$ and three quadrics $q_1, q_2, q_3$ such that 
\[
	\begin{pmatrix}f \\ g \\ h\end{pmatrix}
	=
	\begin{pmatrix}
		q_1 & l_1 \\
		q_2 & l_2 \\
		q_3 & l_3 \\
	\end{pmatrix}
	\begin{pmatrix}x \\ y^2\end{pmatrix}.
\]
If $\codim{q_1,q_2,q_3} \leq 2$, then we apply Lemma~\ref{lemma:3-quadrics}. 
Either the quadrics $q_1,q_2,q_3$ can be expressed in terms of 8 linear forms, or two of the quadrics share a common factor, say $q_1=uz$ and $q_2=vz$ with linear forms $u,v,z$. 
In the former case we have $\pdfgh \leq 13$. 
Namely, $f,g,h$ are in the subring generated by $x,y,l_1,l_2,l_3$ and the 8 linear forms needed to express $q_1,q_2,q_3$. 

In the latter case we are left with eight linear forms $x,y,l_1,l_2,l_3,u,v,z$ and one quadric $q_3$. 
If $q_3$ is in the ideal generated by these eight linear forms, then it can be expressed in terms of those using another set of eight linear forms. 
So $f,g,h$ are in a subring generated by at most 16 linear forms and $\pdfgh \leq 16$. 
And if $q_3 \notin (x,y,l_1,l_2,l_3,u,v,z)$, then $q_3$ is a non-zerodivisor modulo this ideal, that is, the generators of $(x,y,l_1,l_2,l_3,u,v,z)$ together with $q_3$ form a regular sequence of length at most 9 and therefore $\pdfgh \leq 9$. 

Lastly, we need to consider the case $\codim{q_1,q_2,q_3} = 3$ where $q_1,q_2,q_3$ form a regular sequence. 
If they also do so modulo the ideal $(x,y,l_1,l_2,l_3)$, then we have $\pdfgh \leq 8$, as the generators of $(x,y,l_1,l_2,l_3)$ along with $q_1,q_2,q_3$ form a regular sequence of length at most 8.
So we may assume that the images $\bar q_1,\bar q_2,\bar q_3 \in R/(x,y,l_1,l_2,l_3)$ generate an ideal of height $\leq 2$. 
Note that each $\bar{q_i}$ can be lifted back to $q_i$ using 5 linear forms $w_{i1},\ldots,w_{i5}$. 

By Lemma~\ref{lemma:3-quadrics}, either the quadrics $\bar{q_1},\bar{q_2},\bar{q_3}$ can be expressed in terms of 8 linear forms, or two of them share a common factor, say $\bar{q_1}=uz$ and $\bar{q_2}=vz$ with linear forms $u,v,z$. 
In the former case we can place $f,g,h$ in a subring generated by 28 linear forms: 8 linear forms used to express $\bar{q_1},\bar{q_2},\bar{q_3}$, along with $x,y,l_1,l_2,l_3$ and $w_{ij}$ with $i=1,2,3$ and $j=1 \ldots 5$. 
Thus, $\pdfgh \leq 28$. 

In the latter case we have $q_1,q_2 \in k[x,y,l_1,l_2,l_3,u,v,z,w_{1j},w_{2j}]$ with $j=1 \ldots 5$.
Consequently, $f$ and $g$ are contained in this subring as well.
To obtain $h$, we need to adjoin $q_3$.
If $q_3$ is not in the ideal $(x,y,l_1,l_2,l_3,u,v,z,w_{1j},w_{2j})$, then the generators of this ideals along with $q_3$ form a regular sequence of length at most 19 and $\pdfgh \leq 19$.
And if $q_3$ is in the ideal generated by these 18 linear forms, then it can be expressed in terms of those using another set of 18 linear forms. 
Thus, $\pdfgh \leq 36$.
\end{pf}

With Theorem~\ref{thm:unmix-quad} and Proposition~\ref{prop:bound-36}, we are now able to bound the projective dimension of $R/\fgh$ by 36 in the case of multiplicity four. 
By the associativity formula~\eqref{eqn:assoc} there are eleven possible types for the unmixed part $I$, namely: 
\begin{align*}
	\type{4}{1}, 		&& \type{1}{4}, \\
	\type{1,3}{1,1},		&& \type{1,1}{1,3}, \\
	\type{2,2}{1,1},		&& \type{1,1}{2,2}, \\
	\type{1,1,2}{1,1,1},	&& \type{1,1,1}{1,1,2}, \\
	\type{2}{2},			&& \type{1,2}{2,1}, \\
	\type{1,1,1,1}{1,1,1,1}.
\end{align*}
By virtue of Proposition~\ref{prop:bound-36} we may dismiss five of these; we know that $\pdfgh \leq 36$ whenever the length of $R/I$ is at least two locally at an associated prime of multiplicity one. 
There are five such cases which are listed in the right column above. 
In what follows we consider the remaining six cases.

\medskip%

\noindent
\underline{\type{4}{1}}~ 
If $I$ contains a quadric, then $\pdfgh \leq 4$ by Theorem~\ref{thm:unmix-quad}. 
So suppose $I$ does not contain any quadrics; 
in particular, $I$ is non-degenerate. 
By Theorem~\ref{thm:brodmann} of Brodmann and Schenzel, $I$ is the defining ideal of a generic projection of the Veronese surface $V_5 \subset \p5$ onto $\p4$ and it is generated by seven cubics (in 5 variables). 
As $f,g,h$ are linear combinations of those cubics, we have $\pdfgh \leq 5$. 

\medskip%

\noindent
\underline{\type{1,3}{1,1}}~ 
$I=(x,y) \cap P$ with independent linear forms $x,y$ and a height two prime ideal $P$ of multiplicity three. 
If $P$ contains a linear form $l$, then $I$ contains a quadric --- such as $xl$ or $yl$ --- and $\pdfgh \leq 4$ by Theorem~\ref{thm:unmix-quad}. 
If on the other hand $P$ is non-degenerate, then it is the ideal of $2 \times 2$ minors of a $3 \times 2$ matrix of indeterminates, that is, $P$ is generated by three quadrics in at most six variables. 
As $\fgh \subseteq I \subset P$, the three cubics $f,g,h$ can be expressed in terms of those quadrics using no more than nine linear coefficients. 
Thus, $\pdfgh \leq 15$.

\medskip%

\noindent
\underline{\type{2,2}{1,1}}~ 
$I$ is the intersection $(l_1,q_1) \cap (l_2,q_2)$ of two prime ideals where $l_1, l_2$ are linear forms and $q_1, q_2$ are irreducible quadrics. 
As the quadric $l_1 l_2$ belongs to $I$, we have $\pdfgh \leq 4$ by Theorem~\ref{thm:unmix-quad}.

\medskip%

\noindent
\underline{\type{1,1,2}{1,1,1}}~ 
$I$ is the intersection $(x,y) \cap (u,v) \cap (l,q)$ of three prime ideals where $q$ is an irreducible quadric and $x,y,u,v,l$ are (not necessarily independent) linear forms. 
If $\codim{x,y,u,v} = 3$, then, without loss of generality, we may replace $u$ by $x$ and write $I = (x,yv) \cap (l,q)$. 
In this case $I$ contains the quadric $xl$ and $\pdfgh \leq 4$ by Theorem~\ref{thm:unmix-quad}. 

If on the other hand $\codim{x,y,u,v} = 4$, then $I \subset (xu,xv,yu,yv)$ and the cubics $f,g,h$ can be expressed in terms of the quadrics $xu,xv,yu,yv$ using no more than 12 linear forms. 
Thus, $\pdfgh \leq 16$. 

\medskip%

\noindent
\underline{\type{2}{2}}~ 
$I$ is primary to a prime ideal $P = (l,q)$ with a linear form $l$ and an irreducible quadric $q$ such that $\lambda\big( R_{_P}/I_{_P} \big)=2$. 
Thus, locally at $P$, we must have $P_{\!_P}^2 \subset I_{_P}$. 
But primary ideals are contracted ideals in the sense that $I = IR_{_P} \cap R$. 
Hence $P^2 \subset I$ globally.
So $I$ contains the quadric $l^2$ and we have $\pdfgh \leq 4$ by Theorem~\ref{thm:unmix-quad}.

\medskip%

\noindent
\underline{\type{1,1,1,1}{1,1,1,1}}~ 
$I$ is the intersection of four height two prime ideals, each of which is generated by two linear forms. 
So the generators of $I$ are expressed entirely in terms of at most eight (not necessarily independent) linear forms. 
If $I$ contains a quadric, then $\pdfgh \leq 4$ by Theorem~\ref{thm:unmix-quad}. 
And if $I$ is generated in degrees 3 and higher, then the cubics $f,g,h$ are linear combinations (with field coefficients) of the cubic generators of $I$, in which case $\pdfgh \leq 8$.

\subsection{Multiplicity five}
\label{subsection:mult5}

We call to mind the following theorem which is similar in nature to Theorem~\ref{thm:unmix-quad}. 

\begin{thm}[Engheta {\cite[Theorem~17]{e}}]		\label{thm:link-quad}
Let $R$ be a polynomial ring over a field and let $J \subset R$ be an ideal generated by three cubics with $e(R/J) \leq 5$. 
Denote by $I$ the unmixed part of $J$ and let $I'$ be an ideal which is linked to $I$ via cubics. 
If $I'$ contains a quadric, then $\pd{R/J} \leq 4$. 
\end{thm}

Before proceeding with the case of multiplicity five, we single out the following argument which we will employ multiple times in this section as well as in the next. 
Note that there is no assumption on the multiplicity of the ideal. 

\begin{remark}	\label{remark:contained-in-12}
Let $Q$ be an ideal primary to $(x,y)$ with independent linear forms $x,y$ and let $p_1,\ldots,p_k$ be cubic forms in $Q$. 
Suppose $Q \subseteq (x,y)^2+(ax+by)$ with elements $a,b \in \m$ such that $x,y,a,b$ form a regular sequence. 
(In particular, $\deg(ax+by) \geq 2$.) 
Then either the cubics $p_1,\ldots,p_k$ can be expressed entirely in terms of $4(k+1)$ linear forms, or $Q$ is of the form $(x,y)^{e(R/Q)}+(a'x+b'y)$ with quadrics $a',b'$ such that $x,y,a',b'$ form a regular sequence and $\pd{R/Q} \leq 3$. 
\end{remark}
\begin{pf}
The proof of the claim is mainly based on the inclusion 
\[
	(p_1,\ldots,p_k) ~ \subseteq ~ Q ~ \subseteq ~ (x,y)^2+(ax+by).
\]
The only obstacle occurs when $\deg(ax+by)=3$, in which case $a$ and $b$ are quadrics and may involve an arbitrary large number of linear forms. 

Suppose $\deg(ax+by)=3$. 
We first consider the case where one of the $p_i$ has a non-zero contribution from the term $ax+by$, that is, if we write 
\begin{equation}		\label{eqn:pi}
	p_i ~ = ~ l_{i1} \, x^2 + l_{i2} \, xy + l_{i3} \, y^2 + \alpha_i \, (ax+by), 
	\qquad%
	i=1,\ldots,k
\end{equation}
with linear forms $l_{ij}$ and scalars $\alpha_i \in k$, then $\alpha_i$ is non-zero for some $i$. 
Say $\alpha_1 \neq 0$. 
In this case we write $p_1$ as 
\begin{align}	\label{eqn:p1=a'x+b'y}
	p_1 & ~ = ~ l_{11} \, x^2 + l_{12} \, xy + l_{13} \, y^2 + \alpha_1 \, (ax+by) 
	\nonumber \\
		& ~ = ~ \underbrace{(\alpha_1 a + l_{11} x)}_{a'} \, x + 
		\underbrace{(\alpha_1 b + l_{12} x + l_{13} y)}_{b'} \, y, 
\end{align}
and we note that since the elements $x,y,a,b$ form a regular sequence and $\alpha_1 \neq 0$, the elements $x,y,a',b'$ form a regular sequence as well. 
By \cite[Lemma~10]{e} the ideal $(x,y)^{e(R/Q)} + (a'x+b'y)$ is unmixed of multiplicity $e(R/Q)$ and by \cite[Lemma~8]{e} it is equal to $Q$. 
By \cite[Lemma~10]{e} we also have $\pd{R/Q} \leq 3$. 

If on the other hand $\alpha_i = 0$ for all $i=1 \ldots k$, then $(p_1,\ldots,p_k) \subset (x,y)^2$ and by \eqref{eqn:pi} the cubics $p_i$ can be expressed entirely in terms of $3k+2$ linear forms $l_{ij},x,y$. 
Note that the same holds when $\deg(ax+by) \geq 4$. 
We also find ourselves in a similar situation when $\deg(ax+by) = 2$. 
Namely, the cubics $p_i$ are then contained in an ideal generated by four quadrics $x^2,xy,y^2,ax+by$ and so they can be expressed entirely in terms of $4k+4$ linear forms $l_{i1},l_{i2},l_{i3},l_{i4},x,y,a,b$ with $i=1 \ldots k$. 
\end{pf}

We now establish a bound of 20 for the projective dimension of $R/\fgh$ in the case of multiplicity five. 
Let $p_1,p_2$ be any two cubics in the unmixed part $I$ of $\fgh$ which form a regular sequence and let $I'$ denote the link $(p_1,p_2):I$. 
We have $e(R/I')=9-5=4$. 
By the associativity formula~\eqref{eqn:assoc} there are eleven possible types for the link $I'$, namely: 
\begin{align*}
	\type{4}{1}, 		&& \type{1}{4}, \\
	\type{1,3}{1,1},		&& \type{1,1}{1,3}, \\
	\type{2,2}{1,1},		&& \type{1,1}{2,2}, \\
	\type{1,1,2}{1,1,1},	&& \type{1,1,1}{1,1,2}, \\
	\type{2}{2},			&& \type{1,2}{2,1}, \\
	\type{1,1,1,1}{1,1,1,1}.
\end{align*}
The argument which we are about to enter consists of the following parts: 
\begin{itemize}
	\item[$\bullet$] Either the link $I'$ contains a quadric, in which case $\pdfgh \leq 4$ by Theorem~\ref{thm:link-quad}. 
	\item[$\bullet$] Or we give a bound for the projective dimension of $R/I'$ which in turn bounds (by one more) the projective dimension of $R/\fgh$. 
	\item[$\bullet$] Or, by drawing on Remark~\ref{remark:contained-in-12} or by exhibiting that $I'$ is contained in an ideal generated by a set of given quadrics, we show that the cubics $p_1$ and $p_2$ can be expressed entirely in terms of at most 12 linear forms, whereas any one cubic in $I'$ requires at most 8 linear forms. 
\end{itemize}

Recall that $p_1$ and $p_2$ are two arbitrary cubics in $I'$ which form a regular sequence. 
So, unless we are able to obtain a bound for $\pdfgh$ from the first two parts of the above argument, we apply the third part to the choice of, say, $f,g$ and then to $h$ and thus place the cubics $f,g,h$ inside a subring generated by no more than 12+8 linear forms. 
Hence $\pdfgh \leq 20$. 

\medskip%

\noindent
\underline{\type{4}{1}}~ 
If $I'$ contains a quadric, then $\pdfgh \leq 4$ by Theorem~\ref{thm:link-quad}. 
So suppose $I'$ does not contain any quadrics. 
By Theorem~\ref{thm:brodmann} of Brodmann and Schenzel, $I'$ is the defining ideal of a generic projection of the Veronese surface $V_5 \subset \p5$ and $\pd{R/I'} = 4$. 
Thus, $\pdfgh \leq 5$ by Theorem~\ref{thm:general}. 

We point out that the bound of 5 obtained in this case is in fact sharp. 
We will demonstrate this by constructing an example in Section~\ref{section:pdfgh=5}. 

\medskip%

\noindent
\underline{\type{1}{4}}~ 
$I'$ is primary to $(x,y)$ with independent linear forms $x,y$ such that $\lambda\!\left( (R/I')_{(x,y)} \right) = 4$. 
So the Hilbert function of $(R/I')_{(x,y)}$ is either $(1,2,1)$ or $(1,1,1,1)$. 

First suppose $(R/I')_{(x,y)}$ has Hilbert function $(1,2,1)$. 
Then the Hilbert function of $\left( \sfrac{R}{I':(x,y)} \right)_{(x,y)}$ is either $(1,1)$ or $(1,2)$, depending on whether or not $(R/I')_{(x,y)}$ has a socle element outside $(x,y)^2_{(x,y)}$. 

If $\left( \sfrac{R}{I':(x,y)} \right)_{(x,y)}$ has Hilbert function $(1,2)$, then $I':(x,y)=(x,y)^2$ and since $I' \subset I':(x,y) = (x^2,xy,y^2)$, the cubics $p_1,p_2 \in I'$ can be expressed entirely in terms of 8 linear forms. 

If on the other hand $\left( \sfrac{R}{I':(x,y)} \right)_{(x,y)}$ has Hilbert function $(1,1)$, then by Proposition~\ref{prop:class} we have $I':(x,y)=(x,y)^2+(ax+by)$ with elements $a,b$ such that $\codim{x,y,a,b}>3$. 
If the term $ax+by$ is linear, then $I'$ contains quadrics --- such as $(ax+by)x$ and $(ax+by)y$ --- and $\pdfgh \leq 4$ by Theorem~\ref{thm:link-quad}. 
And if $\deg(ax+by) \geq 2$, then we are done by Remark~\ref{remark:contained-in-12}. 

Now suppose $(R/I')_{(x,y)}$ has Hilbert function $(1,1,1,1)$. 
Then the Hilbert function of $\left( \sfrac{R}{I':(x,y)^2} \right)_{(x,y)}$ is $(1,1)$ and by Proposition~\ref{prop:class} we have $I':(x,y)^2=(x,y)^2+(ax+by)$ with elements $a,b$ such that $\codim{x,y,a,b}>3$. 
Again, if $\deg(ax+by) \geq 2$, then we are done by Remark~\ref{remark:contained-in-12}. 

If $\deg(ax+by)=1$, then we may relabel the term $ax+by$ as $x$ so that $I':(x,y)^2=(x,y^2)$. 
In particular, $x(x,y)^2 = (x^3,x^2y,xy^2) \subset I'$. 
Since $(R/I')_{(x,y)}$ has Hilbert function $(1,1,1,1)$, $I'$ must also contain a generator of the form $cx+dy$ with $(c,d) \not\subset (x,y)$. 
Multiplying $cx+dy$ with $y^2$ and reducing it modulo $xy^2$, we see that $dy^3 \in I'$. 
As $(R/I')_{(x,y)}$ has Hilbert function $(1,1,1,1)$, we cannot have $(x,y)^3 \subseteq I'$. 
But $I'$ already contains $(x^3,x^2y,xy^2)$. 
So $y^3 \notin I'$ and therefore $d \in (x,y)$. 
(Recall that $I'$ is primary to $(x,y)$.) 
In particular, $dxy \in (x^2y,xy^2) \subset I'$. 
Multiplying $cx+dy$ with $x$ and reducing it modulo $dxy$, we see that $cx^2 \in I'$. 
As $(c,d) \not\subset (x,y)$ and $d \in (x,y)$, we have $c \notin (x,y)$ and so $x^2 \in I'$. 
Thus, $I'$ contains a quadric and $\pdfgh \leq 4$ by Theorem~\ref{thm:link-quad}. 

\medskip%

\noindent
\underline{\type{1,3}{1,1}}~ 
$I'=(x,y) \cap P$ with independent linear forms $x,y$ and a height two prime ideal $P$ of multiplicity three. 
If $P$ contains a linear form $l$, then $I'$ contains a quadric --- such as $xl$ or $yl$ --- and $\pdfgh \leq 4$ by Theorem~\ref{thm:link-quad}. 
If on the other hand $P$ is non-degenerate, then it is the ideal of $2 \times 2$ minors of a $3 \times 2$ matrix of indeterminates, that is, $P$ is generated by three quadrics in at most six variables. 
As $I' \subset P$, the cubics $p_1,p_2 \in I'$ can be expressed entirely in terms of 12 linear forms. 

\medskip%

\noindent
\underline{\type{1,1}{1,3}}~ 
$I' = (u,v) \cap I_3$ with independent linear forms $u,v$ and an ideal $I_3$ of type $\type{1}{3}$. 
That is, $I_3$ is primary to $(x,y)$ with independent linear forms $x,y$ and $\lambda\!\left( (R/I_3)_{(x,y)} \right) = 3$. 
In particular, $(x,y)^3 \subset I_3$ and the Hilbert function of $(R/I_3)_{(x,y)}$ is either $(1,2)$ or $(1,1,1)$. 
We know that $\codim{x,y,u,v} \geq 3$. 
If $\codim{x,y,u,v} = 4$, then $I' \subset (u,v) \cap (x,y) = (xu,xv,yu,yv)$ and the cubics $p_1,p_2 \in I'$ can be expressed entirely in terms of 12 linear forms. 
So we may assume $\codim{x,y,u,v} = 3$ and without loss of generality, we may replace $u$ by $x$ and write $I' = (x,v) \cap I_3$. 

If $(R/I_3)_{(x,y)}$ has Hilbert function $(1,2)$, then $I_3 = (x,y)^2$ and $I'$ equals $(x^2,xy,y^2v)$. 
It is easily seen that $R/I'$ is \cm{}. 
Consequently, $\pd{R/I'}=2$ and we have $\pdfgh \leq 3$ by Theorem~\ref{thm:general}. 

If on the other hand $(R/I_3)_{(x,y)}$ has Hilbert function $(1,1,1)$, then the quotient $I_3:(x,y)$ is of type \type{1}{2}. 
By Proposition~\ref{prop:class} we have $I_3:(x,y) = (x,y)^2+(ax+by)$ with elements $a,b$ such that $\codim{x,y,a,b}>3$. 

If $\deg(ax+by)=1$, then $I_3=(x^2, xy, y^3, cx+dy^2)$ by \cite[Lemma~13]{e}. 
In particular, modulo $(x,v)$ the ideal $I_3$ is generated by two elements: $(x,v)+I_3 = (x,v)+(y^3,dy^2)$. 
To bound the projective dimension of $R/I'$, we consider the short exact sequence 
\begin{equation}		\label{eqn:ses-xvI3}
	0 \longrightarrow \frac{R}{I'} \longrightarrow \underbrace{\frac{R}{(x,v)} \oplus \frac{R}{I_3}}_{\mathrm{proj.\;dim.}\,\leq\,3} \longrightarrow \underbrace{\frac{R}{(x,v,y^3,dy^2)}}_{\mathrm{proj.\;dim.}\,\leq\,4} \longrightarrow 0
\end{equation}
and note that by \cite[Lemma~12]{e} the middle term has projective dimension $\leq 3$, while the right term is easily seen to have projective dimension $\leq 4$. 
It follows from \eqref{eqn:ses-xvI3} that $\pd{R/I'} \leq 3$, and so $\pdfgh \leq 4$ by Theorem~\ref{thm:general}. 

If $\deg(ax+by) \geq 2$, then we apply the argument of Remark~\ref{remark:contained-in-12} to the ideal $I_3$. 
That is, unless the cubics $p_1,p_2 \in I' \subset I_3$ can be expressed entirely in terms of 12 linear forms, we have $I_3 = (x,y)^3 + (a'x+b'y)$. 
As above, we observe that modulo $(x,v)$ the ideal $I_3$ is generated by two elements: $(x,v)+I_3 = (x,v)+(y^3,b'y)$. 
So we have a short exact sequence similar to \eqref{eqn:ses-xvI3} 
\[
	0 \longrightarrow \frac{R}{I'} \longrightarrow \frac{R}{(x,v)} \oplus \frac{R}{I_3} \longrightarrow \frac{R}{(x,v,y^3,b'y)} \longrightarrow 0
\]
in which the middle term has projective dimension $\leq 3$ by \cite[Lemma~10]{e}, and the right term is easily seen to have projective dimension $\leq 4$. 
As above, $\pd{R/I'} \leq 3$ and $\pdfgh \leq 4$. 

\medskip%

\noindent
\underline{\type{2,2}{1,1}}~ 
$I' = (l_1,q_1) \cap (l_2,q_2)$ with linear forms $l_1,l_2$ and irreducible quadrics $q_1,q_2$. 
As $I'$ contains the quadric $l_1 l_2$, we have $\pdfgh \leq 4$ by Theorem~\ref{thm:link-quad}. 

\medskip%

\noindent
\underline{\type{1,1}{2,2}}~ 
By Proposition~\ref{prop:class} we have 
\[
	I' = (x^2,xy,y^2,ax+by) \cap (u^2,uv,v^2,cu+dv)
\] 
where $x,y,u,v$ are linear forms and $\codim{x,y,u,v} = 3$ or $4$. 
If $\codim{x,y,u,v} = 3$, then, without loss of generality, we may replace $u$ by $x$. 
In this case $I'$ contains the quadric $x^2$ and $\pdfgh \leq 4$ by Theorem~\ref{thm:link-quad}. 
If on the other hand $\codim{x,y,u,v} = 4$, then $I' \subset (x,y) \cap (u,v) = (xu,xv,yu,yv)$. 
So the cubics $p_1,p_2 \in I'$ can be expressed entirely in terms of 12 linear forms. 

\medskip%

\noindent
\underline{\type{1,1,2}{1,1,1}}~ 
$I' = (x,y) \cap (u,v) \cap (l,q)$ with linear forms $x,y,u,v,l$ and an irreducible quadric $q$. 
If $\codim{x,y,u,v} = 3$, then, without loss of generality, we may replace $u$ by $x$ and write $I' = (x,yv) \cap (l,q)$. 
In this case $I'$ contains the quadric $xl$ and $\pdfgh \leq 4$ by Theorem~\ref{thm:link-quad}. 
If on the other hand $\codim{x,y,u,v} = 4$, then $I' \subset (x,y) \cap (u,v) = (xu,xv,yu,yv)$ and the cubics $p_1,p_2 \in I'$ can be expressed entirely in terms of 12 linear forms. 

\medskip%

\noindent
\underline{\type{1,1,1}{1,1,2}}~ 
By Proposition~\ref{prop:class}, $I'$ admits a primary decomposition of the form $I' = (u,v) \cap (s,t) \cap (x^2,xy,y^2,ax+by)$ with linear forms $u,v,s,t,x,y$. 
If $\codim{u,v,s,t}=4$, then $I' \subset (u,v) \cap (s,t) = (us,ut,vs,vt)$ and the cubics $p_1,p_2 \in I'$ can be expressed entirely in terms of 12 linear forms. 

If on the other hand $\codim{u,v,s,t}=3$, then, without loss of generality, $u=s$ and $I' = (u,vt) \cap (x^2,xy,y^2,ax+by)$. 
Note that if $u \in (x,y)$, then $I'$ contains the quadric $u^2$ and $\pdfgh \leq 4$ by Theorem~\ref{thm:link-quad}. 
So we may further assume that $\codim{u,x,y}=3$. 
We now use the inclusion $I' \subset (u,vt) \cap (x,y)$ to bound the number of linear forms needed to write $p_1$ and $p_2$. 

If $vt \notin (x,y)$, then $I' \subset (ux,uy,vtx,vty)$ and the cubics $p_1,p_2 \in I'$ can be expressed entirely in terms of 9 linear forms. 
If on the other hand $vt \in (x,y)$, then either $v \in (x,y)$ or $t \in (x,y)$, for $(x,y)$ is a prime ideal. 
Say $v \in (x,y)$ and, without loss of generality, relabel $v$ as $x$. 
Now $I' \subset (ux,uy,xt)$ and $p_1,p_2 \in I'$ can be expressed entirely in terms of 10 linear forms. 

\medskip%

\noindent
\underline{\type{2}{2}}~ 
$I'$ is primary to a prime ideal $P = (l,q)$ with a linear form $l$ and an irreducible quadric $q$ such that $\lambda\big( R^{}_{_P}/I'_{_P} \big)=2$. 
Thus, locally at $P$, we must have $P_{\!_P}^2 \subset I'_{_P}$. 
But primary ideals are contracted ideals in the sense that $I' = I'R_{_P} \cap R$. 
Hence $P^2 \subset I'$ globally.
So $I'$ contains the quadric $l^2$ and therefore $\pdfgh \leq 4$ by Theorem~\ref{thm:link-quad}.

\medskip%

\noindent
\underline{\type{1,2}{2,1}}~ 
By Proposition~\ref{prop:class}, $I'$ admits a primary decomposition of the form $I' = (x^2,xy,y^2,ax+by) \cap (l,q)$ with linear forms $x,y,l$, an irreducible quadric $q$, and elements $a,b$ such that $\codim{x,y,a,b}>3$. 
If $l \in (x,y)$ or if $\deg(ax+by)=1$, then $I'$ contains the quadric $l^2$ or $(ax+by)l$, respectively, and $\pdfgh \leq 4$ by Theorem~\ref{thm:link-quad}. 
So we may assume that $\codim{x,y,l}=3$ and $\deg(ax+by) \geq 2$, that is, $x,y,l$ and $x,y,a,b$ are both regular sequences. 

As laid out in the proof of Remark~\ref{remark:contained-in-12}, we may further reduce to the case where $\deg(ax+by)=3$ and $ax+by=p_1$. 
(Recall that $I'$ is linked to $I=\fgh\unm$ via two cubics $p_1$ and $p_2$, that is, $I'=(p_1,p_2):I$.) 
Indeed, if $\deg(ax+by)=2$ or $\geq 4$, then the cubics $p_1,p_2 \in I'$ can be expressed entirely in terms of (at most) 12 linear forms. 
The same holds when $\deg(ax+by)=3$ as long as $(p_1,p_2) \subset (x,y)^2$. 
And if $\deg(ax+by)=3$ and one of the cubics, say $p_1$, has a non-zero contribution from the term $ax+by$, then we may replace $ax+by$ by $p_1$ without changing the ideal $(x,y)^2 + (ax+by)$ --- cf.~\eqref{eqn:p1=a'x+b'y} et seq.~on page~\pageref{eqn:p1=a'x+b'y}. 
So without loss of generality $ax+by=p_1$. 

Having replaced the cubic $ax+by$ by $p_1$, we may no longer assume that $a$ and $b$ are reduced modulo $(x,y)$. 
However, as $p_1 \in I'$, we now have $ax+by \in (l,q)$, say $ax+by=cl+l'q$ with a quadric $c$ and a linear form $l'$. 
This reduces the challenge of having to deal with three quadrics $a$, $b$, $q$ to that of having to deal with only two quadrics $c$ and $q$. 
By \cite[Lemma~15]{e} we have 
\[
	I' \;=\; \left[ (x,y)^2 \cap (l,q) \right] + (cl+l'q) ~ \subset ~ (x,y) \cap (l,q). 
\]

To bound the projective dimension of $R/I'$, first suppose $q \in (x,y)$, say $q=l_1x+l_2y$ with linear forms $l_1,l_2$. 
Since $cl+l'q \in (x,y)$, it follows that $cl \in (x,y)$ and as $x,y,l$ form a regular sequence, we must have $c \in (x,y)$, say $c=l_3x+l_4y$ with linear forms $l_3,l_4$. 
Now we can place the generators of $I'$ inside the subring $k[x,y,l,l',l_1,l_2,l_3,l_4]$. 
So $\pd{R/I'} \leq 8$ and $\pdfgh \leq 9$ by Theorem~\ref{thm:general}. 

Now suppose $q \notin (x,y)$. 
Since we may reduce $q$ modulo $l$ without changing the ideal $(l,q)$, this is tantamount to $q \notin (x,y,l)$, that is, $x,y,l,q$ form a regular sequence. 
Thus, from $ax+by=cl+l'q$ we glean $c \in (x,y,q)$, say $c=l_1x+l_2y+\alpha q$ with linear forms $l_1,l_2$ and a scalar $\alpha \in k$. 
This places the generators of $I'$ inside the subring $k[x,y,l,l',l_1,l_2,q]$. 
Let $L$ denote the ideal generated by the linear forms $x,y,l,l',l_1,l_2$. 

If $q \notin L$, then the generators of $L$ along with $q$ form a regular sequence of length at most 7, in which case $\pd{R/I'} \leq 7$ and $\pdfgh \leq 8$. 
If on the other hand $q \in L$, then $q$ can be expressed in terms of the generators of $L$ using no more than six additional linear forms, in which case $\pd{R/I'} \leq 12$ and $\pdfgh \leq 13$. 

\medskip%

\noindent
\underline{\type{1,1,1,1}{1,1,1,1}}~ 
$I'$ is the intersection of four height two prime ideals, each of which is generated by two linear forms. 
Clearly, $\pd{R/I'} \leq 8$ and by Theorem~\ref{thm:general} we have $\pdfgh \leq 9$.

\subsection{Multiplicity six}
\label{subsection:mult6}

Using linkage and Theorem~\ref{thm:general} as our main tools, we give a bound of 20 for the projective dimension of $R/\fgh$ in the case of multiplicity six. 
Let $p_1,p_2$ be any two cubics in the unmixed part $I$ of $\fgh$ which form a regular sequence and let $I'$ denote the link $(p_1,p_2):I$. 
We have $e(R/I')=9-6=3$. 
By the associativity formula~\eqref{eqn:assoc} there are five possible types for the link $I'$, namely: 
\begin{align*}
	\type{3}{1},			&& \type{1}{3}, \\
	\type{1,2}{1,1},		&& \type{1,1}{1,2}, \\
	\type{1,1,1}{1,1,1}.
\end{align*}
In what follows we consider each of these cases and either exhibit a bound for the projective dimension of $R/I'$, and thereupon for that of $R/\fgh$, or we infer that the cubics $f,g,h$ are contained in an ideal generated by a known number of quadrics which are expressed in terms of a fixed number of linear forms.

\medskip%

\noindent
\underline{\type{3}{1}}~ 
$I'$ is a height two prime of multiplicity three. 
Thus, $R/I'$ is \cm{} with $\pd{R/I'}=2$, and $\pdfgh \leq 3$ by Theorem~\ref{thm:general}. 

\medskip%

\noindent
\underline{\type{1}{3}}~ 
$I'$ is primary to $(x,y)$, where $x,y$ are independent linear forms, and $\lambda\!\left( (R/I')_{(x,y)} \right) = 3$. 
Either $I'=(x,y)^2$ or, locally at $(x,y)$, the Hilbert function of $(R/I')_{(x,y)}$ is given by $(1,1,1)$. 
In the former case $R/I'$ is \cm{} and we have $\pdfgh \leq 3$ by Theorem~\ref{thm:general}. 
In the latter case Proposition~\ref{prop:class} yields that $I':(x,y) = (x,y)^2+(ax+by)$ with elements $a,b$ such that $\codim{x,y,a,b}>3$. 
Recall that $I' = (p_1,p_2):I$. 
Thus, we have the following inclusion for any two cubics $p_1,p_2$ in the unmixed part $I$ of $\fgh$ which form a regular sequence: 
\[
	(p_1,p_2) ~ \subset ~ I' ~ \subset ~ I':(x,y) ~ = ~ (x,y)^2+(ax+by).
\]
(Here the elements $x,y,a,b$ depend on the choice of the cubics $p_1$ and $p_2$.) 
We give a bound for $\pdfgh$ by considering the degree of the term $ax+by$. 

If $\deg(ax+by)=1$ for some choice of $p_1$ and $p_2$, then, by \cite[Lemma~13]{e}, $I' = (x^2,xy,y^3,cx+dy^2)$ with elements $c$ and $d$ such that $\codim{x,y,c,d}>3$. 
Thus, $\pd{R/I'} \leq 3$ by \cite[Lemma~12]{e} and $\pdfgh \leq 4$ by Theorem~\ref{thm:general}. 

If $\deg(ax+by) \geq 2$ for some choice of $p_1$ and $p_2$, then we are in the position to invoke an argument which was already used in Section~\ref{subsection:mult5}. 
By Remark~\ref{remark:contained-in-12}, either $\pd{R/I'} \leq 3$ and consequently $\pdfgh \leq 4$, or the cubics $p_1,p_2$ can be expressed in terms of 12 linear forms. 
So, unless $\pdfgh \leq 4$, \emph{every} pair of cubics $p_1,p_2 \in I$ which form a regular sequence can be expressed entirely in terms of 12 linear forms, while any single cubic in $I$ can be expressed entirely in terms of 8 linear forms. 
Thus, $f,g,h$ can be written entirely in terms of 20 linear forms and $\pdfgh \leq 20$. 

\medskip%

\noindent
\underline{\type{1,2}{1,1}}~ 
$I' = (x,y) \cap (l,q)$ with linear forms $x,y,l$ and an irreducible quadric $q$. 
It was shown in \cite[Section~4,~Case~3]{e} that either $\codim{x,y,l,q}=3$ and $R/I'$ is \cm{}, or $\codim{x,y,l,q}=4$ and $\pd{R/I'}=3$. 
Hence $\pdfgh \leq 4$. 

\medskip%

\noindent
\underline{\type{1,1}{1,2}}~ 
By Proposition~\ref{prop:class}, $I'$ admits a primary decomposition of the form $(u,v) \cap (x^2,xy,y^2,ax+by)$ with independent linear forms $u,v$, independent linear forms $x,y$, and elements $a,b$ such that $\codim{x,y,a,b}>3$. 
As so often, we study this intersection through the short exact sequence 
\begin{equation}		\label{eqn:ses-deg6}
	0 \to \frac{R}{I'} \to 
	\underbrace{\frac{R}{(u,v)} \oplus \frac{R}{(x,y)^2+(ax+by)}}_{\mathrm{projective\;dimension}\;\leq\;3} 
	\to \underbrace{\frac{R}{(u,v)+(x,y)^2+(ax+by)}}_{\mathrm{projective\;dimension}\;\leq\;5} \to 0
\end{equation}
in which the projective dimension of the middle term is $\leq 3$ by \cite[Lemma~10]{e}, and the projective dimension of the right term is easily verified to be $\leq 5$. 
(The right term has projective dimension 5 unless either $\codim{u,v,x,y}=3$, or $(a,b) \subset (u,v,x,y)$, or $\deg(ax+by)=1$.) 
Thus, $\pd{R/I'} \leq 4$ by \eqref{eqn:ses-deg6} and $\pdfgh \leq 5$ by Theorem~\ref{thm:general}. 

\medskip%

\noindent
\underline{\type{1,1,1}{1,1,1}}~ 
$I'$ is the intersection of three height two prime ideals, each of which is generated by two linear forms. 
Clearly, $\pd{R/I'} \leq 6$ and by Theorem~\ref{thm:general} we have $\pdfgh \leq 7$.

\section{Three cubics of projective dimension 5}	\label{section:pdfgh=5}

In this section we construct an ideal generated by three cubic forms whose projective dimension equals 5. 
While this answers the question whether an ideal generated by three cubic forms can have projective dimension greater than 4, it is not known whether this is the largest value possible. 

Our construction, which was motivated by part~(c) of the following theorem, leads to an ideal of multiplicity five and corresponds to the case in Section \ref{subsection:mult5} where the link $I'$ of the unmixed part $I$ is of type \type{4}{1}. 
Note that an upper bound of 5 was established in that particular case. 

\begin{thm}[Brodmann, Schenzel {\cite[Theorem 2.1]{bs}}]	\label{thm:brodmann}
A non-degenerate, irreducible projective variety $V$ of multiplicity 4 and codimension 2, which is not a cone, is one of the following: 
\begin{enumerate}
	\item[(a)] a complete intersection of two quadric hypersurfaces 
	\item[(b)] $\dim V \leq 4$ with Betti diagram 
\[ \begin{array}{l|rrrr}
	  & 0 & 1 & 2 & 3 \\
	\hline
	0 & 	1 & - & - & - \\
	1 & - & 1 & - & - \\
	2 & - & 3 & 4 & 1 \\
\end{array} \]
	\item[(c)] (The exceptional case) 
	$V$ is a generic projection of the Veronese surface $V_5 \subset \p5$ with Betti diagram 
\[ \begin{array}{l|rrrrr}
	  & 0 & 1 & 2 & 3 & 4 \\
	\hline
	0 & 	1 & - &  - & - & - \\
	1 & - & - &  - & - & - \\
	2 & - & 7 & 10 & 5 & 1 \\
\end{array} \]
\end{enumerate}
\end{thm}

The starting point of our construction is $I(V_5)$, the defining ideal of the Veronese surface $V_5 \subset \p5$. 
Note that $\codim{I(V_5)}=3$. 
In order to obtain an ideal of height two, we project $V_5$ from a general point of $\p5$ onto $\p4$ and denote the defining ideal of the resulting variety by $I'$. 
(This notation is consistent with that of Section \ref{subsection:mult5}, as $I'$ will be linked to the unmixed part of the three cubics we are about to construct.) 
By part~(c) of Theorem~\ref{thm:brodmann}, $I'$ is generated by seven cubics and $\pd{R/I'}=4$. 
Now, if $I'$ is linked to the unmixed part $I$ of an ideal generated by three cubic forms $f,g,h$, then it follows from Theorem~\ref{thm:general} that $\pdfgh = \pd{R/I'}+1 = 5$. 

To construct an ideal $I$ which is linked to $I'$, we choose two generic cubics $p_1,p_2 \in I'$ and set $I := (p_1,p_2):I'$. 
In the computation carried out below using the computational algebra program Macaulay 2~\cite{m2}, the resulting ideal $I$ is generated by five cubics. 
Choosing $f,g,h$ as three generic linear combinations of these five cubics yields an ideal with $\fgh\unm = I$ and hence $\pdfgh=5$. 

\begin{diagram}[size=2.4em]
\begin{matrix}
	\\
	I(V_5) \ \leadsto \ \\
	~
\end{matrix} &
\begin{matrix}
	e=4 \\
	I' \\
	~
\end{matrix} &
&
&
&
\begin{matrix}
	e=5 \\
	I \\
	~
\end{matrix} &
\begin{matrix}
	\\
	\ = \ \fgh\unm \\
	~
\end{matrix} \\
&
&
\rdLine &
&
\ruLine &
&
& \\
&
&
&
\begin{matrix}
	(p_1,p_2) \\
	e=9
\end{matrix}
\end{diagram}

\bigskip%

\scriptsize \renewcommand{\baselinestretch}{0.9} \begin{verbatim}
Macaulay 2, version 0.9.95
with packages: Classic, Core, Elimination, IntegralClosure,
               LLLBases, Parsing, PrimaryDecomposition,
               SchurRings, TangentCone

i1 : S = QQ[y_0..y_5];

i2 : veronese = trim minors(2, genericSymmetricMatrix(S, y_0, 3))

             2                                    2                       2
o2 = ideal (y  - y y , y y  - y y , y y  - y y , y  - y y , y y  - y y , y  - y y )
             4    3 5   2 4    1 5   2 3    1 4   2    0 5   1 2    0 4   1    0 3

o2 : Ideal of S

i3 : Sbar = S/veronese;

i4 : R = QQ[x_0..x_4];

i5 : link = trim kernel map(Sbar, R, random(Sbar^{1}, Sbar^5));

o5 : Ideal of R

i6 : degrees link

o6 = {{3}, {3}, {3}, {3}, {3}, {3}, {3}}

o6 : List

i7 : p1p2 = ideal(mingens link * random(R^7, R^2));

o7 : Ideal of R

i8 : unmix = p1p2 : link;

o8 : Ideal of R

i9 : degrees unmix

o9 = {{3}, {3}, {3}, {3}, {3}}

o9 : List

i10 : fgh = ideal(mingens unmix * random(R^5, R^3));

o10 : Ideal of R

i11 : top fgh == unmix

o11 = true

i12 : betti res fgh

             0 1 2  3 4 5
o12 = total: 1 3 8 10 5 1
          0: 1 . .  . . .
          1: . . .  . . .
          2: . 3 .  . . .
          3: . . .  . . .
          4: . . 8 10 5 1

o12 : BettiTally
\end{verbatim} \renewcommand{\baselinestretch}{1.1} \normalsize

\bigskip%

Certain outputs of the above computation --- in particular, the output of the cubics $f,g,h$ in line {\tt o10} --- were purposefully suppressed, for the generic choice of the coefficients renders a printout of the resulting polynomials infeasible. 
Yet, to provide the reader with a somewhat manageable example, we repeat the above computation, this time over the finite field $\mathbb{Z}_3 = \mathbb{Z}/3\mathbb{Z}$ rather than the rationals $\mathbb{Q}$, and obtain the following example. 

\begin{example*}
Let $R = \mathbb{Z}_3[X_0,\ldots,X_4]$ and consider the cubic forms 
\[ \begin{split}
f =\ & X_0^3-X_0^2X_2+X_0X_1X_2+X_0X_2^2+X_1X_2^2-X_0^2X_3-X_1X_2X_3-X_2^2X_3-\\
     & X_1X_3^2-X_2X_3^2-X_3^3+X_0^2X_4-X_0X_1X_4-X_1^2X_4-X_1X_2X_4-\\
     & X_0X_3X_4+X_1X_3X_4+X_2X_3X_4-X_3^2X_4+X_2X_4^2+X_3X_4^2+X_4^3, 
\\
g =\ & X_0X_1^2-X_1^3+X_0^2X_2-X_0X_1X_2-X_0X_2^2-X_1X_2^2+X_0^2X_3-\\
     & X_0X_1X_3-X_1X_2X_3+X_2^2X_3+X_0X_3^2-X_3^3-X_0X_1X_4-\\
     & X_1^2X_4+X_0X_2X_4+X_1X_2X_4+X_0X_3X_4+X_2X_3X_4+X_3^2X_4, 
\\
h =\ & X_0^2X_1-X_1^3-X_0^2X_2+X_0X_1X_2-X_1^2X_2-X_0^2X_3-X_1X_2X_3+X_1X_3^2+\\
     & X_0X_1X_4+X_1^2X_4+X_0X_2X_4-X_1X_2X_4+X_0X_3X_4+X_1X_3X_4+X_1X_4^2. 
\end{split} \]
Then $R/\fgh$ has Betti diagram as in line {\tt o12} above. 
In particular, the projective dimension of $R/\fgh$ equals 5. 
\end{example*}

As a caveat, it is worth noting that when performing the above computation over the finite field $\mathbb{Z}_p$, one should verify that the ideal {\tt link} generated in line {\tt o5} --- which is the defining ideal of the projection of the Veronese surface from a general point of $\p5$ onto $\p4$ --- is indeed generic, that is, it is generated by seven cubic forms. 
This check is performed in line {\tt i6}.

\begin{ack}
Thanks are due to Craig Huneke, Hal Schenck, Mike Stillman, and Irena Peeva. 
The computations for the preparation of this paper were performed using the computer algebra system Macaulay 2~\cite{m2}. 
\end{ack}


\begin{thebibliography}{EH}

\bibitem[B]{b}
Burch, L., 1968. 
A note on the homology of ideals generated by three elements in local rings. 
Proc.~Cambridge Philos.~Soc.~64, 949--952.

\bibitem[BS]{bs}
Brodmann, M., Schenzel, P., 2006. 
On varieties of almost minimal degree in small codimension. 
J.~Algebra 305, 789--801.

\bibitem[EH]{eh}
Eisenbud, D., Harris, J., 1987. 
On varieties of minimal degree (a centennial account). 
Proc.~Sympos.~Pure Math.~46, 3--13.

\bibitem[E1]{thesis}
Engheta, B., 2005. 
Bounds on projective dimension. 
Ph.D.~thesis, University of Kansas, Lawrence, KS. 

\bibitem[E2]{e}
Engheta, B., 2007. 
On the projective dimension and the unmixed part of three cubics. 
J.~Algebra 316, 715--734.

\bibitem[M2]{m2}
Grayson, D.R., Stillman, M.E. 
Macaulay 2, a software system for research in algebraic geometry. 
http://www.math.uiuc.edu/Macaulay2/

\bibitem[PS]{ps}
Peeva, I., Stillman, M., 2009.
Open problems on syzygies and Hilbert functions.
J.~Commut.~Algebra 1, no.~1, 159--195.

\end{thebibliography}
\end{document}